\documentclass[12pt]{birkjour}
\usepackage{amsmath,amssymb,amsbsy,amsfonts,amsthm,latexsym,amsopn,amstext,
                            amsxtra,euscript,amscd,tikz}

\newfont{\teneufm}{eufm10}
\newfont{\seveneufm}{eufm7}
\newfont{\fiveeufm}{eufm5}
\newfam\eufmfam
                 \textfont\eufmfam=\teneufm \scriptfont\eufmfam=\seveneufm
                 \scriptscriptfont\eufmfam=\fiveeufm
\def\bbbc{{\mathchoice {\setbox0=\hbox{$\displaystyle\rm C$}\hbox{\hbox
to0pt{\kern0.4\wd0\vrule height0.9\ht0\hss}\box0}}
{\setbox0=\hbox{$\textstyle\rm C$}\hbox{\hbox
to0pt{\kern0.4\wd0\vrule height0.9\ht0\hss}\box0}}
{\setbox0=\hbox{$\scriptstyle\rm C$}\hbox{\hbox
to0pt{\kern0.4\wd0\vrule height0.9\ht0\hss}\box0}}
{\setbox0=\hbox{$\scriptscriptstyle\rm C$}\hbox{\hbox
to0pt{\kern0.4\wd0\vrule height0.9\ht0\hss}\box0}}}}
\def\bbbq{{\mathchoice {\setbox0=\hbox{$\displaystyle\rm
Q$}\hbox{\raise
0.15\ht0\hbox to0pt{\kern0.4\wd0\vrule height0.8\ht0\hss}\box0}}
{\setbox0=\hbox{$\textstyle\rm Q$}\hbox{\raise
0.15\ht0\hbox to0pt{\kern0.4\wd0\vrule height0.8\ht0\hss}\box0}}
{\setbox0=\hbox{$\scriptstyle\rm Q$}\hbox{\raise
0.15\ht0\hbox to0pt{\kern0.4\wd0\vrule height0.7\ht0\hss}\box0}}
{\setbox0=\hbox{$\scriptscriptstyle\rm Q$}\hbox{\raise
0.15\ht0\hbox to0pt{\kern0.4\wd0\vrule height0.7\ht0\hss}\box0}}}}
\def\bbbt{{\mathchoice {\setbox0=\hbox{$\displaystyle\rm
T$}\hbox{\hbox to0pt{\kern0.3\wd0\vrule height0.9\ht0\hss}\box0}}
{\setbox0=\hbox{$\textstyle\rm T$}\hbox{\hbox
to0pt{\kern0.3\wd0\vrule height0.9\ht0\hss}\box0}}
{\setbox0=\hbox{$\scriptstyle\rm T$}\hbox{\hbox
to0pt{\kern0.3\wd0\vrule height0.9\ht0\hss}\box0}}
{\setbox0=\hbox{$\scriptscriptstyle\rm T$}\hbox{\hbox
to0pt{\kern0.3\wd0\vrule height0.9\ht0\hss}\box0}}}}
\def\bbbs{{\mathchoice
{\setbox0=\hbox{$\displaystyle     \rm S$}\hbox{\raise0.5\ht0\hbox
to0pt{\kern0.35\wd0\vrule height0.45\ht0\hss}\hbox
to0pt{\kern0.55\wd0\vrule height0.5\ht0\hss}\box0}}
{\setbox0=\hbox{$\textstyle        \rm S$}\hbox{\raise0.5\ht0\hbox
to0pt{\kern0.35\wd0\vrule height0.45\ht0\hss}\hbox
to0pt{\kern0.55\wd0\vrule height0.5\ht0\hss}\box0}}
{\setbox0=\hbox{$\scriptstyle      \rm S$}\hbox{\raise0.5\ht0\hbox
to0pt{\kern0.35\wd0\vrule height0.45\ht0\hss}\raise0.05\ht0\hbox
to0pt{\kern0.5\wd0\vrule height0.45\ht0\hss}\box0}}
{\setbox0=\hbox{$\scriptscriptstyle\rm S$}\hbox{\raise0.5\ht0\hbox
to0pt{\kern0.4\wd0\vrule height0.45\ht0\hss}\raise0.05\ht0\hbox
to0pt{\kern0.55\wd0\vrule height0.45\ht0\hss}\box0}}}}
\def\bbbz{{\mathchoice {\hbox{$\sf\textstyle Z\kern-0.4em Z$}}
{\hbox{$\sf\textstyle Z\kern-0.4em Z$}}
{\hbox{$\sf\scriptstyle Z\kern-0.3em Z$}}
{\hbox{$\sf\scriptscriptstyle Z\kern-0.2em Z$}}}}

\newtheorem{theorem}{Theorem}[section]

\newtheorem{corollary}[theorem]{Corollary}
\def\cB{{\mathcal B}}
\def\cC{{\mathcal C}}

\def\cU{{\mathcal U}}

\def\({\left(}
\def\){\right)}
\def\[{\left[}
\def\]{\right]}
\def\<{\langle}
\def\>{\rangle}

\def\F{\mathbb{F}}

\def\vec#1{\mathbf{#1}}
\def\dist{\mathrm{dist}}

\def\fC{\mathfrak{C}}

\begin{document}

\title{Solutions to polynomial congruences in well-shaped sets}
\author{Bryce Kerr}
\address{Department of Computing, Macquarie University, Sydney, NSW 2109, Australia}

\email{bryce.kerr@mq.edu.au}

\begin{abstract}
We use a generalization of Vinogradov's mean value theorem of S. Parsell, S. Prendiville and T. Wooley  and ideas of W. Schmidt to give nontrivial bounds for the number of solutions to polynomial congruences, when the solutions lie in a very general class of sets, including all convex sets.
\end{abstract}

\maketitle
\section{Introduction}
Given integer $m$ and a polynomial $F(X_1,\dots,X_d)\in \mathbb{Z}_m[X_1,\dots X_d]$  and some $\Omega \subseteq [0,1]^{d}$, we let $N_F(\Omega)$ denote the number of solutions $\mathbf{x}=(x_1,\dots x_d)\in \mathbb{Z}^d$ to the congruence
\begin{equation}
\label{solutions}
F(\mathbf{x})\equiv 0 \text{ \ (mod $m$)}
 \ \text{ \ with \ } \  
\left(\frac{x_1}{m},\dots \frac{x_d}{m}\right)\in \Omega.
\end{equation}
 Questions concerning the distribution of solutions to polynomial congruences have been considered in a number of works (for example \cite{CCGHSZ,GrShZa,Shp1,Zum}). In \cite{Fouv} Fouvry gives an asymptotic formula for the number of solutions to systems of polynomial congruences in small cubic boxes for a wide class of systems (see also \cite{FoKa,Luo,ShpSk,Skor}). Shparlinski \cite{Sp} uses the results of \cite{Fouv} and ideas of \cite{Schm}  to obtain an asymptotic formula for the number of solutions to the same systems when the solutions lie in a very general class of sets. For the case of a single polynomial, $F$ in $d$ variables, Shparlinski \cite{Sp} shows that for suitable $\Omega$, when the modulus $m=p$ is prime,
\begin{equation}
\label{shp}
N_F(\Omega)=p^{d-1}(\mu(\Omega)+O(p^{-1/4}\log{p}))
\end{equation}
provided $F$ is irreducible over $\mathbb{C}$ and $\mu$ denotes the Lebesgue measure on $[0,1]^d$. This gives an asymptotic formula for $N_F(\Omega)$ provided $\mu(\Omega)\geq p^{-1/4+\epsilon}$ and a nontrivial upper bound for $N_F(\Omega)$ when $\mu(\Omega)\ge p^{-5/4+\varepsilon}$. We follow the method of \cite{Sp} to give an upper bound for $N_F(\Omega)$ without any restrictions on our polynomial $F$ when the modulus $m$ is composite.  We first establish an upper bound for $N_F(\Omega)$  when $\Omega$ is a  cube. This gives a generalization of Theorem $1$ of \cite{cp}. Although we follow the same argument, the difference is our use of a multidimensional version of Vinogradov's mean value theorem (Theorem 1.1 of \cite{mdws}). To extend the bound from cubes to more general sets $\Omega$, we approximate $\Omega$  by cubes using ideas based on Theorem $2$ of \cite{Schm}. \\
\section{Definitions}
\indent We let $\mu$ denote  the Lebesgue measure on $[0,1]^d$,  $||.||$ the Euclidian norm and define the distance between  $\mathbf{x}\in [0,1]^d$ and $\Omega\subseteq [0,1]^d$ to be
$$\dist(\mathbf{x},\Omega)=\displaystyle\inf_{\mathbf{y}\in \Omega}||\mathbf{x}-\mathbf{y}||.$$
As in \cite{Sp}, we say that  $\Omega \subseteq [0,1]^{d}$ is \emph{well-shaped} if there exists $C=C(\Omega)$ such that for every $\varepsilon>0$ the measures of the sets 
$$
\Omega_\varepsilon^{+} = \left\{ \vec{u} \in  [0,1]^d \backslash
\Omega \ : \ \dist(\vec{u},\Omega) < \varepsilon \right\},
$$
$$
\Omega_\varepsilon^{-} = \left\{ \vec{u} \in \Omega \ : \
\dist(\vec{u},[0,1]^d \backslash \Omega )  < \varepsilon  \right\}
$$
exist and satisfy
\begin{equation}
\label{well-shaped}
\mu(\Omega_{\varepsilon}^{\pm})\leq C\varepsilon.
\end{equation}
 From Lemma 1 of \cite{Schm} all convex subsets of $[0,1]^d$ are well-shaped and from equation $(2)$ of  \cite{Weyl}, if the boundary of $\Omega$ is a manifold of dimension $n-1$ with bounded surface area  then $\Omega$ is well-shaped, for suitably chosen $C.$ \\ 
\indent For $\mathbf{x}=(x_1,\dots x_d)$ we write $a\leq \mathbf{x} \leq b$ if  $a \leq x_1, \dots, x_d \leq b$.  Given a $d$-tuple of non-negative integers
 $\mathbf{i}=(i_1,i_2,\dots i_d)$, we set $\mathbf{x}^{\mathbf{i}}=x_{1}^{i_1}x_{2}^{i_2}\dots x_{d}^{i_d}$
and $|\mathbf{i}|=i_1+i_2+\dots i_d$. We let $r$ denote the number of distinct $d$-tuples, $\mathbf{i}$ with
$1\leq |\mathbf{i}|\leq k,$ so that
\begin{equation}
\label{r} 
r=\binom{k+d}{d}-1.
\end{equation}
\indent We will always suppose $m$ is an integer greater than $2$. Given $F\in \mathbb{Z}_m[X_1,\dots,X_d]$, we let $k$ denote the degree of $F$ and $d$ the number of variables. Writing
$$F(\mathbf{x})=\displaystyle\sum_{0\leq |\mathbf{i}| \leq k}\beta_{\mathbf{i}}\mathbf{x}^{\mathbf{i}}, \quad \beta_{\mathbf{i}}\in \mathbb{Z}_m$$
we define $$g_F=\displaystyle\min_{|\mathbf{i}|=k}\gcd(m,\beta_i).$$
\indent We use $g(t) \ll f(t)$ and $g(t)=O(f(t))$ to mean that there exists some absolute constant $\alpha$ such that $|g(t)|\leq \alpha f(t)$ for all values of $t$ within some specified range. Whenever we use $\ll$ and $O$, unless stated otherwise  the implied constant  will depend only on $d$, $k$ and the particular $C$ in $(\ref{well-shaped}).$ Similarily $o(1)$ denotes a term which is sufficiently small when our parameter is large enough in terms of $d$, $k$ and $C$. 

\section{Main Results}
We can now present our main results:
\begin{theorem}
\label{multi}
For positive $K_1, \dots, K_d, L,H,R\ge 1$, integer $m$ and 
$$F(\mathbf{x})=\displaystyle\sum_{0\leq |\mathbf{i}| \leq k}\beta_{\mathbf{i}}\mathbf{x}^{\mathbf{i}}\in \mathbb{Z}_m[X_1,\dots, X_d ]$$ of degree $k\ge2$ with $g_F=1$, let $M_F(H,R)$ denote the number of solutions to  the congruence
\begin{equation}
\label{polynomial equation 1}
F(\mathbf{x})\equiv y \ \   (\mathrm{mod} \  m)
\end{equation}
 with
$$(\mathbf{x},y)\in [K_1+1,K_1+H]\times \dots \times [K_d+1,K_d+H]\times [L+1,L+R].$$
Then uniformly over all $K_1,\dots,K_d,L\ge 1$
$$M_F(H,R)\le H^d\left(\left(\frac{R}{H^k}\right)^{1/2r(k+1)}+\left(\frac{R}{m}\right)^{1/2r(k+1)}\right)m^{o(1)}$$
as $H\rightarrow \infty.$
\end{theorem}
Arguing from heuristics, we expect the bound for $M_F(H,R)$ to be around
$$M_F(H,R)\le H^d\left(\frac{R}{m}\right)m^{o(1)}$$
which can be directly compared with Theorem~1. Similarly, by considering the first term in Theorem~1  we immediatley see when this bound for $M_F(H,R)$ is worse than the trivial bound $M_F(H,R)\le H^d$. \newline \indent
 Also, if $m=p$ is prime and $F[X_1,\dots, X_d]$ is not multilinear, i.e $F$ is not linear in each of its variables, then Theorem~\ref{multi} is trivial. This may be seen by the following argument. First we may show by slightly adjusting the proof of Theorem~1 of~\cite{cp} that for $G \in \mathbb{Z}_p[X]$ of degree $k\ge 2$
\begin{equation}
\label{aabb}
M_G(H,R)\le H\left(\left(\frac{R}{H^k}\right)^{1/2k(k+1)}+\left(\frac{R}{p}\right)^{1/2k(k+1)}\right)p^{o(1)}.
\end{equation}
 Supposing $F\in \mathbb{Z}_p[X_1,\dots,X_{d}]$ is not multilinear, then after re-ordering the variables we may suppose for some $k_0 \ge 2$ that
\begin{equation}
\label{induction step}
F[X_1,\dots,X_{d}]=\sum_{i=0}^{k_0}X_{d}^iF_{i}[X_1,\dots,X_{d-1}]
\end{equation}
with $F_{k_0}\neq 0$ and consider separately the values of $X_1,\dots, X_{d-1}$  such that 
\begin{equation*}
\label{equiv1}
F_{k_0}[X_1,\dots,X_{d-1}]\equiv0 \pmod p
\end{equation*}
 and
\begin{equation*}
\label{notequiv1} 
F_{k_0}[X_1,\dots,X_{d-1}] \not  \equiv 0 \pmod p.
\end{equation*}
For the first case we use the assumption that $p$ is prime and induction on $d$ to bound the number of values $X_1,\dots,X_{d-1}$ such that
$F_{k_0}[X_1,\dots,X_{d-1}]\equiv0 \pmod p$ by $O(H^{d-2})$ and bound the number of solutions to
$$F[X_1,\dots,X_{d}]\equiv y \pmod p$$
in remaining variables $X_{d},Y$ trivially by $RH$. \newline \indent For the second case, we bound the number of $X_1,\dots, X_{d-1}$ such that 
$F_{k_0}[X_1,\dots,X_{d-1}] \not  \equiv 0 \pmod p$ trivially by $H^{d-1}$ and  bound the number of solutions in the remaining variables $X_{d},Y$ by applying~\eqref{aabb} to the equation~\eqref{induction step}. Combining the above two cases gives
$$M_F(H,R)\le  H^d\left(\frac{R}{H}+\left(\frac{R}{H^k}\right)^{1/2k(k+1)}+\left(\frac{R}{m}\right)^{1/2k(k+1)}\right)p^{o(1)}$$
which can be compared with Theorem~1. \newline
\indent Taking $R=1$ in Theorem~\ref{multi} we get,
\begin{corollary}
\label{cubes}
For any cube $B\subseteq [0,1]^{d}$ of side length $\frac{1}{h}$, $F \in \mathbb{Z}_m[X_1,\dots, X_d ]$ of degree $k\ge 2$ with $g_F=1$ we have
$$N_F(B)\le \left(\frac{m}{h}\right)^{d-k/2r(k+1)+o(1)}+m^{d-1/2r(k+1)+o(1)}\left(\frac{1}{h}\right)^{d+o(1)}$$
as \  $\dfrac{m}{h} \rightarrow \infty.$
\end{corollary}
 Taking $R=H$ in Theorem~\ref{multi} we get,
\begin{corollary}
\label{cubes 1}
Suppose $F \in \mathbb{Z}_m[X_1,\dots, X_d ]$ of degree $k\ge 2$ with $g_F=1$ is of the form,
$$F(X_1,\dots, X_d)=G(X_1,\dots X_{d-1})-X_{d}$$
for some $G \in \mathbb{Z}_m[X_1,\dots,X_{d-1}]$, then
for any cube $B\subseteq [0,1]^{d}$ of side length $\frac{1}{h}$, we have
\begin{align*}
N_F(B)&\le \left(\frac{m}{h}\right)^{d-1-(k-1)/2r(k+1)+o(1)}+ m^{d-1+o(1)}\left(\frac{1}{h}\right)^{d-1+1/2r(k+1)+o(1)}
\end{align*}
as \  $\dfrac{m}{h} \rightarrow \infty$, where $r$ corresponds to $d-1$ in the definition~\eqref{r}. 
\end{corollary} 
We use the above Corollaries to estimate $N_F(\Omega)$ for well-shaped $\Omega$.
\begin{theorem}
\label{multi vws}
Suppose $F \in \mathbb{Z}_m[X_1,\dots, X_d ]$ satisfies the conditions of Corollary~\ref{cubes} and $\Omega \subset [0,1]^{d}$ is well-shaped with
$\mu(\Omega)\geq m^{-1}$. Then we have
$$N_F(\Omega)\leq m^{d-k/2r(k+1)+o(1)}\mu(\Omega)^{1-k/2r(k+1)}+m^{d-1/2r(k+1)+o(1)}\mu(\Omega)$$
 as  $m \rightarrow \infty.$ 
\end{theorem}
\begin{theorem}
Suppose $F\in \mathbb{Z}_m[X_1,X_2,\dots X_d]$ satisfies the conditions of Corollary~\ref{cubes 1} and $\Omega \subset [0,1]^{d}$ is well-shaped. Then we have
$$
N_F(\Omega) \leq \begin{cases}m^{d-1+o(1)}\mu(\Omega)^{1/2r(k+1)}, \quad \mu(\Omega)\ge m^{-1+1/k}  \\
m^{d-1-(k-1)/2r(k+1)+o(1)}\mu(\Omega)^{-(k-1)/2r(k+1)}, \quad m^{-1}\le \mu(\Omega)<m^{-1+1/k}.
 \end{cases}
 $$
as $m\rightarrow \infty.$
\end{theorem}

\section{Proof of Theorem 3.1}
Making a change of variables we may assume $(\mathbf{K},L)=(0,\dots,0)$. Suppose for integer $s$ we have $\mathbf{x}_1,\mathbf{x}_2,\dots, \mathbf{x}_{2s}$  satisfying (\ref{polynomial equation 1})  with \\
$\mathbf{x}_j=(x_{j,1},x_{j,2},\dots,x_{j,d}).$ Then
$$ F(\mathbf{x}_1)+ F(\mathbf{x}_2)+\dots+ F(\mathbf{x}_s)- F(\mathbf{x}_{s+1})-\dots - F(\mathbf{x}_{2s})\equiv z \text{\ (mod $m$)} $$
for some \ $-sR \leq z \leq sR.$  Hence there exists  $-sR \leq u \leq sR$ \  such that
\begin{equation}
\label{M bound}
M_F(H,R)^{2s}\leq (1+2sR)T(u,H)
\end{equation}
with $T(u,H)$ equal to the number of solutions to the congruence
\begin{equation}
\label{polynomial equation 2} 
F(\mathbf{x}_1)+ F(\mathbf{x}_2)+\dots+ F(\mathbf{x}_s)- F(\mathbf{x}_{s+1})-\dots - F(\mathbf{x}_{2s})\equiv u \text{\ (mod $m$)}
\end{equation}
with each co-ordinate of $\mathbf{x_j}$ between $1$ and $H.$ \\ \\ Since
$$F(\mathbf{x})=\displaystyle\sum_{0\leq |\mathbf{i}| \leq k}\beta_{\mathbf{i}}\mathbf{x}^{\mathbf{i}}, \text{ \ for some\  $\beta_{\mathbf{i}}\in \mathbb{Z}_m$} $$
we may write (\ref{polynomial equation 2}) in the form
\begin{equation}
\label{linear}
\displaystyle\sum_{1\leq |\mathbf{i}| \leq k}\beta_{\mathbf{i}}\lambda_{\mathbf{i}}\equiv u \text{ \ (mod $m$)}
\end{equation}
with 
\begin{equation}
\label{lambda}
\lambda_{\mathbf{i}}=\mathbf{x}_1^{\mathbf{i}}+\dots+\mathbf{x}_s^{\mathbf{i}}
-\mathbf{x}_{s+1}^{\mathbf{i}}-\dots -\mathbf{x}_{2s}^{\mathbf{i}}.
\end{equation}
Since $g_F=1$, we choose $\mathbf{i}_0$ with $|\mathbf{i}_0|=k$ and $\gcd(\beta_{\mathbf{i}_0},m)=1$. Considering (\ref{linear}) as a linear equation in $\lambda_{\mathbf{i}}$, if we let $\lambda_{\mathbf{i}}$, $\mathbf{i}\neq \mathbf{i}_0$ take arbitrary values then  $\lambda_{\mathbf{i}_0}$ is determined uniquely $\pmod m$. 
 Since we have
\begin{equation}
\label{lambda bound}
-sH^{|\mathbf{i}|}\leq \lambda_{\mathbf{i}}\leq sH^{|\mathbf{i}|}
\end{equation}
there are at most 
\begin{equation}
\label{T bound step 1}
\left(1+(2s+1)\dfrac{H^{k}}{m}\right)\displaystyle\prod_{\substack{\mathbf{i}\neq \mathbf{i}_0\\ 1\leq |\mathbf{i}|\leq k }}(2s+1)H^{|\mathbf{i}|}=(2s+1)^{r-1}H^{K-k}\left(1+(2s+1)\dfrac{gH^{k}}{m}\right)
\end{equation}
solutions to (\ref{linear}) in integer variables $\lambda_{\mathbf{i}},$ with
$$K=\displaystyle\sum_{1\leq |\mathbf{i}| \leq k}|\mathbf{i}|=\frac{d}{d+1}(r+1)k.$$ 
For \ $U=(u_{\mathbf{i}})_{1\leq |\mathbf{i}|\leq k}$ with each $u_{\mathbf{i}}\in \mathbb{Z},$ \ let \ 
$J_{s,k,d}(U, H)$ \ denote the number of solutions in integers, $\lambda_{\mathbf{i}}$, to
\begin{equation}
\label{equation over Z}
\lambda_{\mathbf{i}}=u_{\mathbf{i}}, \text{ \ $1\leq |\mathbf{i}|\leq k$}
\end{equation}
with each $\mathbf{x}_j$ having components between $1$ and $H$ and we write $J_{s,k,d}(U, H)=J_{s,k,d}(H)$ when $U=(0)_{1\leq |\mathbf{i}|\leq k}.$ Let $\cU$ be the collection of sets  $U=(u_{\mathbf{i}})_{1\leq |\mathbf{i}|\leq k}$ such that $|u_{\mathbf{i}}|\leq sH^{|\mathbf{i}|}$ and 
$$\displaystyle\sum_{1\leq |\mathbf{i}| \leq k}\beta_{\mathbf{i}}u_{\mathbf{i}}\equiv u \text{ \ (mod $m$)}$$
so that the cardinality of $\cU$ is bounded by~\eqref{T bound step 1}.
We see that 
\begin{equation}
\label{T bound step 2}
T(u,H)\leq \displaystyle\sum_{U \in \cU}J_{s,k,d}(U,H),
\end{equation}
since if $\mathbf{x}_{0,1} \dots \mathbf{x}_{0,2s}$ is a solution to~\eqref{polynomial equation 2}, then  the integers $\lambda_{0,\mathbf{i}}$, defined by 
$$\lambda_{0,\mathbf{i}}=\mathbf{x}_{0,1}^{\mathbf{i}}+\dots+\mathbf{x}_{0,s}^{\mathbf{i}}
-\mathbf{x}_{0,s+1}^{\mathbf{i}}-\dots -\mathbf{x}_{0,2s}^{\mathbf{i}}, \ \ \  \ \ 1\leq |\mathbf{i}|\leq k$$
are a solution to~\eqref{linear} and the $\mathbf{x}_{0,1} \dots \mathbf{x}_{0,2s}$ are a solution to 
$$\lambda_{\mathbf{i}}=\lambda_{0,\mathbf{i}}, \text{ \ $1\leq |\mathbf{i}| \leq k$}.$$
So if we let $U_0=(\lambda_{0,\mathbf{i}})_{1\leq |\mathbf{i}|\leq k}$, then we see that the solution to \eqref{polynomial equation 2}, $\mathbf{x}_{0,1} \dots \mathbf{x}_{0,2s},$ is counted by the term $J_{s,k,d}(U_0,H)$ in \eqref{T bound step 2}. By~\eqref{T bound step 1} and~\eqref{T bound step 2}, we have 
\begin{equation}
\label{t bound}
T(u,H)\leq (2s+1)^{r-1}H^{K-k}\left(1+(2s+1)\dfrac{H^{k}}{m}\right)J_{s,k,d}(V ,  H)
\end{equation}
for some $V\in \cU$. \ Although for any $U \in \cU$ we have the inequality
$$J_{s,k,d}(U,  H)\leq J_{s,k,d}(H).$$
Since if we let $\boldsymbol{\alpha}=(\alpha_\mathbf{i})_{1\leq |\mathbf{i}| \leq k}$ and
$$S(\boldsymbol{\alpha})=\displaystyle\sum_{1\leq \mathbf{x}\leq H}\exp\left(2\pi i \displaystyle\sum_{1\leq |\mathbf{i}|\leq k}\alpha_\mathbf{i} \mathbf{x}^{\mathbf{i}}\right)$$ 
then for $\lambda_{\mathbf{i}}$ defined as in ~\eqref{lambda} we have 
\begin{align*}
J_{s,k,d}(U,  H)&=\displaystyle\sum_{1\leq \mathbf{x}_1,\dots \mathbf{x}_{2s}\leq H}
\displaystyle\int_{[0,1]^{r}}\exp\left(2\pi i \displaystyle\sum_{1\leq |\mathbf{i}|\leq k}\alpha_{\mathbf{i}} (\lambda_{\mathbf{i}}
-u_{\mathbf{i}})\right)d\boldsymbol{\alpha}
\\ &=
\displaystyle\int_{[0,1]^{r}}|S(\boldsymbol{\alpha})|^{2s}\exp\left(-2\pi i\displaystyle\sum_{1\leq |\mathbf{i}| \leq k}\alpha_{\mathbf{i}}u_{\mathbf{i}}\right)d\boldsymbol{\alpha} \\  
&\leq \displaystyle\int_{[0,1]^{r}}|S(\boldsymbol{\alpha})|^{2s}d\boldsymbol{\alpha}=J_{s,k,d}(H)
\end{align*}
where  the integral is over the variables $\alpha_{\mathbf{i}}$, $1\leq |\mathbf{i}| \leq k$. Hence by ~\eqref{M bound} and ~\eqref{t bound} we have
\begin{equation}
\label{M bound 2}
M_F(H,R)^{2s}\leq  (1+2sR) (2s+1)^{r-1}H^{K-k}\left(1+(2s+1)\dfrac{H^{k}}{m}\right)J_{s,k,d}(  H).
\end{equation}
By Theorem 1.1 of \cite{mdws} we have for  $s\geq r(k+1)$
$$J_{s,k,d}(H) \ll H^{2sd-K+\epsilon}$$
for any $\epsilon>0$ provided $H$ is sufficiently large in terms of $k,d$ and $s$. Inserting this bound into~\eqref{M bound 2} gives
$$M_F(H,R)^{2s}\ll  RH^{K-k}\left(1+\dfrac{H^k}{m}\right)H^{2sd-K+\epsilon}$$
and the result follows taking $s=r(k+1)$. \\
\qed
\section{Proof of Theorem 3.4}
As in \cite{Schm} we begin with choosing $\mathbf{a}=(a_1,\dots a_{d})$ with each co-ordinate irrational. For integer $j$ let $\fC(j)$ be the set of cubes of the form
\begin{equation}
\label{cuubes}
\left[a_1+\frac{u_1}{j},a_1+\frac{u_1+1}{j}\right]\times \dots \times  \left[a_{d}+\frac{u_{d}}{j},a_{d}+\frac{u_{d}+1}{j}\right], \ \ \ u_i \in \mathbb{Z}.
\end{equation}
Since each $a_i$ is irrational, no point (\ref{solutions}) lies in two distinct cubes (\ref{cuubes}). Given integer $M>0$, let $\varepsilon=2d^{\frac{1}{2}}/2^ M$ and consider the set 
$$\Omega_{\varepsilon}=\Omega \cup \Omega_{\varepsilon}^{+}.$$
Since $\Omega$ is well-shaped, we have
\begin{equation}
\label{omega big}
\mu(\Omega_{\epsilon})=\mu(\Omega)+O\left(\frac{1}{2^M}\right).
\end{equation}
Let $\cC(j)$ be the cubes of $\fC(j)$ lying inside $\Omega_\varepsilon$ and we suppose $j\leq 2^M.$ Then by (\ref{omega big}) we obtain,
\begin{equation}
\label{c upper bound}
\# \cC(j)\leq j^d\mu(\Omega_\varepsilon)\leq j^{d}\mu(\Omega)+O\left(\frac{j^d}{2^M}\right)=
j^{d}\mu(\Omega)+O\left(j^{d-1}\right).
\end{equation}
\begin{figure}
\begin{tikzpicture}[scale=0.35]
%\draw [ultra thick] (-8,-7) rectangle (8,9);
\begin{scope}[rotate around={90:(0,0)}]
\draw[ dashed,very thick, rotate around={20:(-0.8,0.8)}] (-0.80,0.8) ellipse (6.77 and 8.8);
\draw[fill=gray!50, very thick, rotate around={20:(-0.8,0.8)}] (-0.80,0.8) ellipse (5.35 and 7.38);

%\cB_1%
\draw  [very thick] (-4,-4) rectangle (4,4);

%\cB_2%

\draw  [very thick] (-4,4) rectangle (0,8);
%\cB_3%

\draw  [very thick] (-4,-4) rectangle (-2,-6);
\draw  [very thick] (2,-4) rectangle (4,-6);
\draw  [very thick] (0,4) rectangle (2,6);
\draw  [very thick](2,4) rectangle (4,6);
\draw  [very thick](0,6) rectangle (2,8);
\draw  [very thick](4,-2) rectangle (6,0);
\draw  [very thick] (0,-4) rectangle (-2,-6);
\draw  [very thick] (2,-4) rectangle (0,-6);
\draw  [very thick](-4,-2) rectangle (-6,0);
\draw  [very thick](-4,0) rectangle (-6,2);
\draw  [very thick] (-4,2) rectangle (-6,4);
\draw  [very thick] (-4,4) rectangle (-6,6);

%\cB_4%
\draw  [very thick] (-1,-6) rectangle (0,-7);
\draw  [very thick] (0,-6) rectangle (1,-7);
\draw  [very thick](-4,6) rectangle (-5,7);
\draw  [very thick] (4,2) rectangle (5,3);
\draw  [very thick] (4,1) rectangle (5,2);
\draw  [very thick] (4,0) rectangle (5,1);
\draw  [very thick](4,-3) rectangle (5,-2);
\draw  [very thick] (4,-4) rectangle (5,-3);
\draw  [very thick] (-5,-4) rectangle (-4,-3);
\draw  [very thick](-5,-3) rectangle (-4,-2);
\draw  [very thick](-6,-3) rectangle (-5,-2);
\draw  [very thick] (-7,0) rectangle (-6,1);
\draw  [very thick](-7,1) rectangle (-6,2);
\draw  [very thick](-7,2) rectangle (-6,3);
\draw  [very thick] (-7,3) rectangle (-6,4);
\draw  [very thick] (1,-6) rectangle (2,-7);
\draw  [very thick] (-6,4) rectangle (-7,5);
\draw  [very thick](-6,-1) rectangle (-7,0);
\draw  [very thick] (-2,-6) rectangle (-1,-7);

\draw [very thick]  (2,6) rectangle (3,7);
\draw  [very thick]  (4,3) rectangle (5,4);
\draw  [very thick] (5,0) rectangle (6,1);
\draw [very thick] (-4,8) rectangle (-3,9);
\draw [very thick] (-5,-5) rectangle (-4,-4);
\draw [very thick] (-6,6) rectangle (-5,7);
\draw [very thick] (-5,7) rectangle (-4,8);
\draw [very thick]  (-7,5) rectangle (-6,6);
\draw [very thick] (-3,8) rectangle (-2,9);
\draw [very thick] (-2,8) rectangle (-1,9);
\draw [very thick] (-1,8) rectangle (0,9);

\draw  [very thick]  (-7,-1) rectangle (-6,-2);

\draw  [very thick] (2,-6) rectangle (3,-7);
\draw  [very thick] (4,-5) rectangle (5,-4);

%coloured%
%\cB_1%
\draw  [very thick,,fill=gray!50,opacity=0.5] (-4,-4) rectangle (4,4);

%\cB_2%

\draw  [very thick,fill=gray!50,opacity=0.5] (-4,4) rectangle (0,8);
%\cB_3%

\draw  [very thick,fill=gray!50,opacity=0.5] (-4,-4) rectangle (-2,-6);
\draw  [very thick,fill=gray!50,opacity=0.5] (2,-4) rectangle (4,-6);
\draw  [very thick,fill=gray!50,opacity=0.5] (0,4) rectangle (2,6);
\draw  [very thick,fill=gray!50,opacity=0.5](2,4) rectangle (4,6);
\draw  [very thick,fill=gray!50,opacity=0.5](0,6) rectangle (2,8);
\draw  [very thick,fill=gray!50,opacity=0.5](4,-2) rectangle (6,0);
\draw  [very thick,fill=gray!50,opacity=0.5] (0,-4) rectangle (-2,-6);
\draw  [very thick,fill=gray!50,opacity=0.5] (2,-4) rectangle (0,-6);
\draw  [very thick,fill=gray!50,opacity=0.5](-4,-2) rectangle (-6,0);
\draw  [very thick,fill=gray!50,opacity=0.5](-4,0) rectangle (-6,2);
\draw  [very thick,fill=gray!50,opacity=0.5] (-4,2) rectangle (-6,4);
\draw  [very thick,fill=gray!50,opacity=0.5] (-4,4) rectangle (-6,6);

%\cB_4%
\draw  [very thick,fill=gray!50,opacity=0.5] (-1,-6) rectangle (0,-7);
\draw  [very thick,fill=gray!50,opacity=0.5] (0,-6) rectangle (1,-7);
\draw  [very thick,fill=gray!50,opacity=0.5](-4,6) rectangle (-5,7);
\draw  [very thick,fill=gray!50,opacity=0.5] (4,2) rectangle (5,3);
\draw  [very thick,fill=gray!50,opacity=0.5] (4,1) rectangle (5,2);
\draw  [very thick,fill=gray!50,opacity=0.5] (4,0) rectangle (5,1);
\draw  [very thick,fill=gray!50,opacity=0.5](4,-3) rectangle (5,-2);
\draw  [very thick,fill=gray!50,opacity=0.5] (4,-4) rectangle (5,-3);
\draw  [very thick,fill=gray!50,opacity=0.5] (-5,-4) rectangle (-4,-3);
\draw  [very thick,fill=gray!50,opacity=0.5](-5,-3) rectangle (-4,-2);
\draw  [very thick,fill=gray!50,opacity=0.5](-6,-3) rectangle (-5,-2);
\draw  [very thick,fill=gray!50,opacity=0.5] (-7,0) rectangle (-6,1);
\draw  [very thick,fill=gray!50,opacity=0.5](-7,1) rectangle (-6,2);
\draw  [very thick,fill=gray!50,opacity=0.5](-7,2) rectangle (-6,3);
\draw  [very thick,fill=gray!50,opacity=0.5] (-7,3) rectangle (-6,4);
\draw  [very thick,fill=gray!50,opacity=0.5] (1,-6) rectangle (2,-7);
\draw  [very thick,fill=gray!50,opacity=0.5] (-6,4) rectangle (-7,5);
\draw  [very thick,fill=gray!50,opacity=0.5](-6,-1) rectangle (-7,0);
\draw  [very thick,fill=gray!50,opacity=0.5] (-2,-6) rectangle (-1,-7);

\draw [very thick,fill=gray!50,opacity=0.5]  (2,6) rectangle (3,7);
\draw  [very thick,fill=gray!50,opacity=0.5]  (4,3) rectangle (5,4);
\draw  [very thick,fill=gray!50,opacity=0.5] (5,0) rectangle (6,1);
\draw [very thick,fill=gray!50,opacity=0.5] (-4,8) rectangle (-3,9);
\draw [very thick,fill=gray!50,opacity=0.5] (-5,-5) rectangle (-4,-4);
\draw [very thick,fill=gray!50,opacity=0.5] (-6,6) rectangle (-5,7);
\draw [very thick,fill=gray!50,opacity=0.5] (-5,7) rectangle (-4,8);
\draw [very thick,fill=gray!50,opacity=0.5]  (-7,5) rectangle (-6,6);
\draw [very thick,fill=gray!50,opacity=0.5] (-3,8) rectangle (-2,9);
\draw [very thick,fill=gray!50,opacity=0.5] (-2,8) rectangle (-1,9);
\draw [very thick,fill=gray!50,opacity=0.5] (-1,8) rectangle (0,9);

\draw  [very thick,fill=gray!50,opacity=0.5]  (-7,-1) rectangle (-6,-2);

\draw  [very thick,fill=gray!50,opacity=0.5] (2,-6) rectangle (3,-7);
\draw  [very thick,fill=gray!50,opacity=0.5] (4,-5) rectangle (5,-4);

\end{scope}
%\node at (-6,5) {$\Omega_{\varepsilon}$};
%\node at (-3,3) {$\Omega$};
%\node at (-3,3) {$\cB_1$};
%\node at (-7,-1) {$\cB_2$};
%\node [above left] at (-4.7,-5.3) {$\cB_3$};
%\node at (0.5,-6.6) {$\cB_4$};
\end{tikzpicture}
\caption{The sets $\Omega_{\varepsilon}$ and $\Omega$  with the corresponding \ $\cB_1,\cB_2,\cB_3,\cB_4$.}
\end{figure}
Also, since a cube of side length $1/j$ has diameter $\varepsilon_j=d^{\frac{1}{2}}/j$, we see that the cubes of $\cC(j)$ cover $\Omega_\epsilon \setminus (\Omega_\varepsilon)_{\varepsilon_j}^{-}$ and hence 
$$ \# \cC(j)\geq j^d\left(\mu(\Omega_\varepsilon)-\mu((\Omega_\varepsilon)_{\varepsilon_j}^{-})\right).$$
But for $j \leq 2^M$, we have 
$$(\Omega_\varepsilon)_{\varepsilon_j}^{-}\subseteq \Omega_{\varepsilon_j}^{-}\cup \Omega_{\varepsilon}^{+}$$
and since $\Omega$ is well-shaped
$$\mu((\Omega_\varepsilon)_{\varepsilon_j}^{-})\leq \mu(\Omega_{\varepsilon_j}^{-})+ \mu(\Omega_{\varepsilon}^{+})
\ll \frac{1}{j}$$
so we get
$$\# \cC(j)\geq j^d\mu(\Omega_\varepsilon)+O(j^{d-1}).$$
Combining this with (\ref{c upper bound}) gives
\begin{equation}
\label{C bound}
\# \cC(j)=j^{d}\mu(\Omega)+O(j^{d-1}) \text{ \  for $j\leq 2^M.$}
\end{equation}
 Let $\cB_1=\cC(2)$ and for $2\leq i \leq M$ we let $\cB_i$ be the set of cubes from $\cC(2^i)$ that are not contained in any cubes from $\cC(2^{i-1}).$ Then we have $\# \cB_1 =\# \cC(2)$ and for $2\leq i \leq M,$ the cubes from both $\cB_i$ and  $\cC(2^{i-1})$ are contained in $\Omega_\varepsilon$. This gives
 \begin{align*}
 \# \cB_i + 2^d\# \cC(2^{i-1})\leq 2^{id}\mu(\Omega)+O\left(\frac{2^{id}}{2^M}\right)\leq  2^{id}\mu(\Omega)+O\left(2^{i(d-1)}\right)
\end{align*}
and hence by (\ref{C bound}) 
\begin{equation}
 \label{B bound}
 \# \cB_i \ll 2^{i(d-1)}.
 \end{equation}
We have
\begin{equation}
\label{omega in}
\Omega  \subseteq \bigcup_{i=1}^M  \bigcup_{\Gamma \in \cB_i}\Gamma
\end{equation}
since if $\mathbf{x} \in \Omega$ then $$\dist(\mathbf{x},[0,1]^d \setminus \Omega_\varepsilon)\geq \varepsilon.$$
Although $\mathbf{x} \in \Gamma$ for some $\Gamma \in \fC(2^M)$ and since $\Gamma$ has diameter $\varepsilon/2$  we have $\Gamma \in \cC(2^M).$ Since the union of the cubes from $\cC(2^{i-1})$ is contained in the union from $\cC(2^{i})$ we get (\ref{omega in}). Hence
$$N_F(\Omega)\leq \displaystyle\sum_{i=1}^{M}\displaystyle\sum_{\Gamma \in \cB_i}N_F(\Gamma)$$
and using Corollary \ref{cubes}, as $m2^{-M}\rightarrow \infty$
\begin{align*}
\displaystyle\sum_{i=1}^{M}\displaystyle\sum_{\Gamma \in \cB_i}N_F(\Gamma) &\ll
\displaystyle\sum_{i=1}^{M}\displaystyle\sum_{\Gamma \in \cB_i}
\left(\frac{m}{2^{i}}\right)^{d-k/2r(k+1)+o(1)} \\ & \ \ \ +\displaystyle\sum_{i=1}^{M}\displaystyle\sum_{\Gamma \in \cB_i}m^{d-1/2r(k+1)+o(1)}2^{-i(d+o(1))} \\
 &\ll m^{d-k/2r(k+1)+o(1)}2^{o(M)}\displaystyle\sum_{i=1}^{M}2^{ik/2r(k+1)} \frac{\# \cB_i }{2^{id}} \\ & \ \ \ + m^{d-1/2r(k+1)+o(1)}2^{o(M)}\displaystyle\sum_{i=1}^{M}\frac{\# \cB_i}{2^{id}}. 
\end{align*}
We use (\ref{omega big}) to bound 
\begin{equation*}
\displaystyle\sum_{i=1}^{M} \frac{\# \cB_i}{2^{id}}\leq \mu(\Omega_{\varepsilon})=\mu(\Omega)+O\left(\frac{1}{2^M}\right)
\end{equation*}
and from (\ref{B bound}),  for $N\leq M$
\begin{align*}
\displaystyle\sum_{i=1}^{M}2^{ik/2r(k+1)} \frac{\# \cB_i }{2^{id}} &= \nonumber 
\displaystyle\sum_{i=1}^{N}2^{ik/2r(k+1)} \frac{\# \cB_i }{2^{id}}+\displaystyle\sum_{i=N+1}^{M}2^{ik/2r(k+1)} \frac{\# B_i }{2^{id}} \\ &\ll 2^{Nk/2r(k+1)}\displaystyle\sum_{i=1}^{N}\frac{\# \cB_i }{2^{id}}
+\displaystyle\sum_{i=N+1}^{M}2^{ik/2r(k+1)}\frac{2^{i(d-1)}}{2^{id}} \\
&\ll 2^{Nk/2r(k+1)}\left(\mu(\Omega)+2^{-M}\right)+2^{-N(1-k/2r(k+1))} \\
&\ll 2^{Nk/2r(k+1)}\left(\mu(\Omega)+2^{-N}\right).
\end{align*}
Hence we get
\begin{align}
\label{optimize}
N_F(\Omega)&\leq m^{d-k/2r(k+1)+o(1)} 2^{Nk/2r(k+1)+o(M)}\left(\mu(\Omega)+2^{-N}\right) \\
& \ \ \ +2^{o(M)}m^{d-1/2r(k+1)+o(1)}\left(\mu(\Omega)+2^{-M}\right). \nonumber
\end{align}
Recalling that $\mu(\Omega)\geq m^{-1},$ to balance the two terms involving $N$, we choose
$$2^{-N}\leq  \mu(\Omega) \log{m} <2^{-N+1}.$$
Substituting this choice into (\ref{optimize}) gives,
\begin{align*}
N_F(\Omega) &\leq m^{d-k/2r(k+1)}2^{o(M)}\mu(\Omega)^{1-k/2r(k+1)} \\
& \ \ \ +m^{d-1/2r(k+1)+o(1)}2^{o(M)}\left(\mu(\Omega)+2^{-M}\right).
\end{align*}
The same choice for $M$ is essentially optimal,
\begin{equation}
\label{M choice}
2^{-M}\leq m^{-1}\log{m}\leq 2^{-M+1}.
\end{equation}
This gives
$$N_F(\Omega)\leq m^{d-k/2r(k+1)+o(1)}\mu(\Omega)^{1-k/2r(k+1)}+m^{d-1/2r(k+1)+o(1)}\mu(\Omega)$$
where we have replaced $2^{o(M)}$ with $m^{o(1)}$ since $\mu(\Omega)\geq m^{-1}.$ Theorem 3.4 follows since for the choice of $M$ in (\ref{M choice}), for $\mu(\Omega)\geq m^{-1}$
$$m2^{-M}\gg m^{-1}\mu(\Omega)\log{m}\geq \log{m}$$  
which tends to infinity as $m\rightarrow \infty.$  \qed \\ \
\section{Proof of Theorem 3.5}
Using the same constructions from Theorem 3.5, we have
$$N_F(\Omega)\leq \displaystyle\sum_{i=1}^{M}\displaystyle\sum_{\Gamma \in \cB_i}N_F(\Gamma).$$
Hence by Corollary~\ref{cubes 1}
\begin{align}
\label{nf}
N_F(\Omega)&\leq 2^{o(M)}m^{d-1-(k-1)/2r(k+1)+o(1)}\displaystyle\sum_{i=1}^{M}2^{i(1+(k-1)/2r(k+1))}\frac{\#\cB_i}{2^{id}} \nonumber \\ & \ \ \ +2^{o(M)}m^{d-1+o(1)}\displaystyle\sum_{i=1}^{M}2^{i(1-1/2r(k+1))}\frac{\#\cB_i}{2^{id}}.
\end{align}
For the first sum by~\eqref{B bound},
\begin{align*}
\displaystyle\sum_{i=1}^{M}2^{i(1+(k-1)/2r(k+1))}\frac{\#\cB_i}{2^{id}}\le \displaystyle\sum_{i=1}^{M}2^{i(k-1)/2r(k+1))}\ll 2^{M(k-1)/2r(k+1)}.
\end{align*}
For the second sum,
\begin{align*}
\displaystyle\sum_{i=1}^{M}2^{i(1-1/2r(k+1))}\frac{\#\cB_i}{2^{id}}&=\displaystyle\sum_{i=1}^{N}2^{i(1-1/2r(k+1))}\frac{\#\cB_i}{2^{id}} +\displaystyle\sum_{i=N+1}^{M}2^{i(1-1/2r(k+1))}\frac{\#\cB_i}{2^{id}} \\
&\ll 2^{N(1-1/2r(k+1))}\left(\mu(\Omega)+\frac{1}{2^M}\right)+2^{-N/2r(k+1)} \\
&\ll 2^{N(1-1/2r(k+1))}\mu(\Omega)+2^{-N/2r(k+1)}.
\end{align*}
Substituting the above bounds into~\eqref{nf} gives
\begin{align*}
N_F(\Omega)&\leq 2^{o(M)}m^{d-1-(k-1)/2r(k+1)+o(1)}2^{M(k-1)/2r(k+1)} \\ & \ \ \ + 2^{o(M)}m^{d-1+o(1)}\left(2^{N(1-1/2r(k+1))}\mu(\Omega)+2^{-N/2r(k+1)}\right).
\end{align*}
For $\mu(\Omega)\ge m^{-1+1/k}$ we choose $N$ to balance the first and last terms then choose $M$ to balance the remaining terms, so that
$$2^{M-1}< \mu(\Omega)^{1/(k-1)}m\le 2^M$$
$$2^{-N}\le 2^{M(k-1)}m^{-(k-1)}<2^{-N+1}$$
which gives $N\le M$ and
$$N_F(\Omega)\le m^{d-1+o(1)}\mu(\Omega)^{1/2r(k+1)}.$$
If $m^{-1}\le \mu(\Omega) <m^{-1+1/k}$ then we choose $N$ to balance the last two terms and take $M$ as small as possible subject to the condition $N\le M$. This gives  $$2^{-M}\le \mu(\Omega)<2^{-M+1}$$ $$N=M$$ 
and
\begin{align*}
N_F(\Omega) &\le m^{d-1-(k-1)/2r(k-1)}\mu(\Omega)^{-(k-1)/2r(k+1)} \\
& \ \ \ +m^{d-1+o(1)}\mu(\Omega)^{1/2r(k+1)}.
\end{align*}
Combining the above two bounds completes the proof. \qed
\section{Comments}
Using the methods of Theorem 3.4 and Theorem 3.5, we have not been able to to give bounds for $N_F(\Omega)$ which are nontrivial when $\mu(\Omega)\leq m^{-1}$. This seems to be caused by two factors, the bound from Corollary \ref{cubes} and the bounds for $\mu(\Omega_{\varepsilon})^{\pm}$, which affect the estimates (\ref{omega big}) and (\ref{B bound}). For certain cases with prime modulus we may be able to do better than Theorem 3.5. For example, the same  method  may be combined with other bounds replacing Corollary \ref{cubes 1} for more specific families of polynomials. This has the potential to obtain sharper estimates for such polynomials and also to increase the range of values of $\mu(\Omega)$ for which an analogue of Theorem 3.5 would apply. For example, Bourgain, Garaev, Konyagin and Shparlinski~\cite{BGKS1} consider the number
$J_{\nu}(p,h,s;\lambda)$ of solutions to the congruence
\begin{equation*}
(x_1+s)\dots(x_{\nu}+s)\equiv \lambda \ \ (\text{mod} \ p), \ \  1\leq x_1, \dots , x_{\nu} \leq h.
\end{equation*}
They show that if $h<p^{1/(\nu^2-1)}$ then we have the bound
\begin{equation}
\label{b1}
J_{\nu}(p,h,s;\lambda)\leq \exp \left(c(\nu)\frac{\log{h}}{\log\log{h}}\right)
\end{equation}
for some constant $c(\nu)$ depending only on $\nu$  (Lemma 2.33 of \cite{BGKS1}). \\
In \cite{BGKS2}, the same authors consider the number
$K_\nu(p,h,s)$ of solutions to the congruence
\begin{equation*}
(x_1+s)\dots(x_\nu+s)\equiv (y_1+s)\dots(y_\nu+s) \not \equiv 0 \ \ (\text{mod} \ p),
\end{equation*}
$$1\leq x_1, \dots, x_\nu, y_1, \dots ,y_\nu \leq h$$
and show that 
\begin{equation}
\label{b2}
K_\nu(p,h,s)\leq \left(\frac{h^{\nu}}{p^{\nu/e_{\nu}}}+1\right)h^{\nu}\exp\left(c(\nu)\frac{\log{h}}{\log\log{h}}\right)
\end{equation}
for some constants $e_{\nu}$ and $c(\nu)$ depending only on $\nu$  (Theorem 17 of \cite{BGKS2}). \\ \indent Another possible way to improve on our results  for certain classes of well-shaped sets is to use Weyl's formula for tubes (equation (2) of  \cite{Weyl}) and Steiner's formula for convex bodies (equation (4.2.27) of \cite{Schn}) to give an explicit constant
in (\ref{well-shaped}) for certain subsets of $[0,1]^d$ for which these formula are valid. This would have the effect of improving on the bounds (\ref{omega big}) and  (\ref{B bound}) and hence the bound for $N_F(\Omega)$ and possibly the range of values of $\mu(\Omega)$ for which this bound would be valid.
\section*{Acknowledgements}
The author would like to thank Igor Shparlinski for suggesting this problem and for his guidance while working on it and writing the current paper.

\end{document}